\documentclass[10pt]{amsart}
\usepackage{amssymb}
\usepackage{amsmath,amssymb,amsfonts,amsthm,
latexsym, amscd, amsfonts, epsfig}
\usepackage{mathrsfs}
\usepackage{color}
\usepackage{eucal}%    caligraphic-euler fonts: \mathcal{ }
\usepackage{eufrak}%   frak-euler        fonts: \mathfrak{ }
\usepackage[all]{xypic}
\usepackage{xspace}

\theoremstyle{plain}
\newtheorem{theorem}{Theorem}

\newtheorem{lemma}{Lemma}

\theoremstyle{definition}
\newtheorem{definition}{Definition}

\theoremstyle{example}

\theoremstyle{remark}

\numberwithin{equation}{section}

\begin{document}

%%%
%%%
%%%%%%%%%%%%%%%%%%%%%%%%%%%%%%%%%%%%%%%%%%%%%%%%%%%%%%%%%%%%%%%%%%%%%%%%%%
%%
%%%
\title[Pseudoknot RNA Structures with Arc-Length $\ge 3$]
      {Pseudoknot RNA Structures with Arc-Length $\ge 3$}
\author{Emma Y. Jin and Christian M. Reidys$^{\,\star}$}
\address{Center for Combinatorics, LPMC-TJKLC \\
         Nankai University  \\
         Tianjin 300071\\
         P.R.~China\\
         Phone: *86-22-2350-6800\\
         Fax:   *86-22-2350-9272}
\email{reidys@nankai.edu.cn}
\thanks{}
\date{June 2007}
\begin{abstract}
In this paper we study $k$-noncrossing RNA structures with
arc-length $\ge 3$, i.e.~RNA molecules in which for any $i$, the
nucleotides labeled $i$ and $i+j$ ($j=1,2$) cannot form a bond and
in which there are at most $k-1$ mutually crossing arcs. Let ${\sf
S}_{k,3}(n)$ denote their number. Based on a novel functional
equation for the generating function $\sum_{n\ge 0}{\sf
S}_{k,3}(n)z^n$, we derive for arbitrary $k\ge 3$ exponential growth
factors and for $k=3$ the subexponential factor. Our main result is
the derivation of the formula ${\sf S}_{3,3}(n) \sim
\frac{6.11170\cdot 4!}{n(n-1)\dots(n-4)} 4.54920^n$.
\end{abstract}
\maketitle
{{\small
%\tableofcontents
}}

%%%
%%%
%%%%%%%%%%%%%%%%%%%%%%%%%%%%%%%%%%%%%%%%%%%%%%%%%%%%%%%%%%%%%%%%%%%%%%%%
%%%
%%%

\section{Introduction}

%%%
%%%%%%%%%%%%%%%%%%%%%%%%%%%%%%%%%%%%%%%%%%%%%%%%%%%%%%%%%%%%%%%%%%%%%%%%
%%%

RNA is believed to be central for the understanding of evolution.
It acts as genotypic legislative in form of viruses and viroids and
as phenotypic executive in form of ribozymes, capable of catalytic
activity, cleaving other RNA molecules.
This dualism gives rise to the hypothesis that RNA may have preceded
DNA and proteins, therefore playing a key role in prebiotic evolution.
In light of growing support of an {\it RNA world} \cite{Schuster:02} and
even RNA-based metabolisms or the prospect of self-replicating RNA 
\cite{Poole:06a} it is the phenotypic aspect of RNA that still lacks 
deeper understanding. 
Despite the fact that pseudoknot RNA structures are known to be of central 
importance \cite{Science:05a}, little is known from a theoretical
point of view, for instance only recently their generating function has been
obtained \cite{Reidys:07pseu}.

Let us provide first some background on RNA sequences and structures which
allows us to put our results into context.
The primary sequence of an RNA molecule is the sequence of
nucleotides {\bf A}, {\bf G}, {\bf U} and {\bf C} together with the
Watson-Crick ({\bf A-U},{\bf G-C}) and ({\bf U-G}) base pairing
rules specifying the pairs of nucleotides can potentially form
bonds. Single stranded RNA molecules form helical structures whose
bonds satisfy the above base pairing rules and which, in many cases,
determine or are even tantamount to their function.
Three decades ago Waterman {\it et.al.} pioneered the concept
of RNA secondary structures
\cite{Penner:93c,Waterman:79a,Waterman:78a,Waterman:94a,Waterman:80},
a concept being best understood when considering structure as a diagram,
drawing the primary sequence of nucleotides
horizontally and ignoring all chemical bonds of its backbone. Then one
draws all bonds, i.e.~nucleotide interactions satisfying the Watson-Crick
base pairing rules (and {\bf G}-{\bf U} pairs) as arcs in the upper
halfplane, thereby identifying structure with the set of all
arcs. Then secondary structures have no two arcs $(i_1,j_1)$, $(i_2,j_2)$, 
where $i_1<j_1$ and $i_2<j_2$ with the property $i_1<i_2<j_1<j_2$ and 
all arcs have at least length $2$.
While the concept of secondary structure is of fundamental
importance it is well-known that there exist additional types of
nucleotide interactions \cite{Science:05a}. These bonds are called
pseudoknots \cite{Westhof:92a} and occur in functional RNA (RNAseP
\cite{Loria:96a}), ribosomal RNA \cite{Konings:95a} and are
conserved in the catalytic core of group I introns.
In plant viral RNAs pseudoknots mimic tRNA structure and in {\it in
vitro} RNA evolution \cite{Tuerk:92} experiments have produced
families of RNA structures with pseudoknot motifs, when binding
HIV-1 reverse transcriptase.

Leaving the paradigm of RNA secondary structures, i.e.~studying RNA structures
with crossing bonds, poses challenging problems for computational biology.
Prediction algorithms for RNA pseudoknot structures are much harder
to derive since there exists no {\it a priori} tree-structure and
the subadditivity of local solutions is not guaranteed. In fact pseudoknot
RNA structures are genuinely non-inductive and seem to be best described by the
mathematical language of vacillating tableaux \cite{Reidys:07pseu,Chen:07a}.
One approach of categorizing RNA pseudoknot structures consists in considering
the maximal size of sets of mutually crossing bonds,
leading to the notion of $k$-noncrossing RNA structures \cite{Reidys:07pseu}.
%%%
%%%%%%%%%%%%%%%%%%%%%%%% Figures ex1 %%%%%%%%%%%%%%%%%%%%%%%%%%%%%%%%%%%%%%
%%%
\begin{figure}[ht]
\centerline{%
\epsfig{file=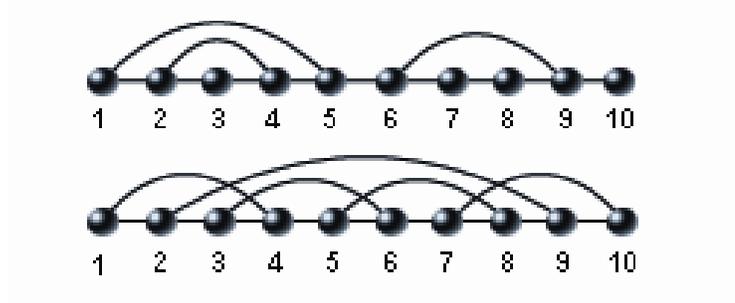,width=0.7\textwidth}\hskip15pt } \caption{\small
$k$-noncrossing RNA structures. (a) secondary structure (with
isolated labels $3,7,8,10$), (b) non-planar $3$-noncrossing
structure } \label{F:4}
\end{figure}
%%%
%%%%%%%%%%%%%%%%%%%%%%%%%%%%%%%%%%%%%%%%%%%%%%%%%%%%%%%%%%%%%%%%%%%%%%%%%
%%%
This concept is very intuitive, a $k$-noncrossing RNA structure has at most
$k-1$ mutually crossing arcs and a minimum bond-length of $2$, i.e.~for any
$i$, the nucleotides $i$ and $i+1$ cannot form a bond.
In this paper we will consider $k$-noncrossing RNA structures with arc-length 
$\ge 3$. Their analysis is based on the generating function, which
has coefficients that are alternating sums. This fact makes even the 
computation of the exponential growth factor a difficult task. 
To make things worse, in case of arc-length $\ge 3$ there exists no 
explicit formula for the coefficients, which could only be computed 
via a recursion formula.

\subsection{Organization and main results}
Let $\mathcal{S}_{k,3}(n)$ ($k\ge 3)$ denote the set of $k$-noncrossing RNA
structures with arc-length $\ge 3$ and let ${\sf S}_{k,3}(n)= \vert
\mathcal{S}_{k,3}(n)\vert $. In Section~\ref{S:pre} we provide the
necessary background on $k$-noncrossing RNA structures with
arc-length $\ge 3$ and their generating function $\sum_{n\ge 0}{\sf
S}_{k,3}(n)z^n$. In Section~\ref{S:exp} we compute the exponential
factor of ${\sf S}_{k,3}(n)$. To make it easily accessible to a broad 
readership we give an elementary proof based on real analysis and 
transformations of the generating function. 
Central to our proof is the functional identity of
Lemma~\ref{L:func} and its generalized version in Lemma~\ref{L:ana}.
In Section~\ref{S:sub} we present the asymptotic analysis of
$\mathcal{S}_{3,3}(n)$, using Flajolet {\it et.al.}'s singular
expansions and transfer theorems
\cite{Flajolet:05,Flajolet:99,Flajolet:94,Popken:53,Odlyzko:92}. 
This analysis is similar to \cite{Reidys:07pseu} but involves solving a
quartic instead of a quadratic polynomial in order to localize the 
singularities. The main result of the paper is 

%%%
%%%%%%%%%%%%%%%%%%%%%%%%%%%%%%%%%%%%%%%%%%%%%%%%%%%%%%%%%%%%%%%%%%%%%%%%%
%%%
{\bf Theorem.}$\,$ {\it The number of $3$-noncrossing RNA structures
with arc length $\ge 3$ is asymptotically given by
\begin{eqnarray}
\label{E:konk3} {\sf S}_{3,3}(n) & \sim & \frac{6.11170\cdot
4!}{n(n-1)\dots(n-4)}\,
4.54920^n \ ,
\end{eqnarray}
where ${\sf s}_{3,3}^{}(n)=\frac{6.11170\cdot 4!}{n(n-1)\dots(n-4)}=
146.6807\left[\frac{1}{ n^5}-\frac{35}{4n^6}+\frac{1525}{32
n^7}+{O}(n^{-8})\right]$.
}
%%%
%%%%%%%%%%%%%%%%%%%%%%%%%%%%%%%%%%%%%%%%%%%%%%%%%%%%%%%%%%%%%%%%%%%%%%%%%
%%%

In the table below we display the quality of our approximation by listing
the subexponential factors, i.e.~we compare for $k=3$ the quantities 
${\sf S}_{k,3}(n)/(4.54920)^n$
obtained from the generating function (Theorem~\ref{T:cool2}), which
are the exact values and the asymptotic expressions ${\sf
s}_{3,3}^{}(n)$, respectively.
\begin{center}
\begin{tabular}{c|c|c|c|c|c}
\hline
  \multicolumn{6}{c}{\textbf{The sub exponential factor}}\\
  \hline
$n$ & ${\sf S}_{3,3}(n)/(4.54920)^n$ & ${\sf s}_{3,3}^{}(n)$ &$n$ & ${\sf S}_{3,3}(n)/(4.54920)^n$ & ${\sf s}_{3,3}^{}(n)$\\
\hline \small 10  & \small $3.016\times 10^{-4}$ & \small$4.851\times 10^{-3}$ &
\small 60  & \small $3.457 \times 10^{-7}$ & \small$2.238\times 10^{-7}$\\
\small 20  & \small $2.017 \times 10^{-5}$ & \small$7.884\times 10^{-5}$ &
\small 70  & \small $1.476\times 10^{-7}$ & \small$1.010\times 10^{-7}$\\
\small 30  & \small $3.513 \times 10^{-6}$ & \small$8.577\times
10^{-6}$ & \small 80  & \small $3.783\times 10^{-8}$ &
\small$5.085\times
10^{-8}$\\
\small 40  & \small $9.646\times 10^{-7}$ & \small$1.858\times
10^{-6}$ & \small 90  & \small $2.154\times 10^{-8}$ & \small$2.781\times 10^{-8}$\\
\small 50  & \small $5.627\times 10^{-7}$ & \small$5.769\times
10^{-7}$ &
\small 100  & \small $1.299\times 10^{-8}$ & \small$1.624\times 10^{-8}$\\
\end{tabular}
\end{center}
%%%
%%%%%%%%%%%%%%%%%%%%%%%%%%%%%%%%%%%%%%%%%%%%%%%%%%%%%%%%%%%%%%%%%%%%%%%%
%%%

\section{$k$-noncrossing RNA structures with arc-length $\ge 3$}\label{S:pre}

\begin{figure}[ht]\label{F:2}
\centerline{ \epsfig{file=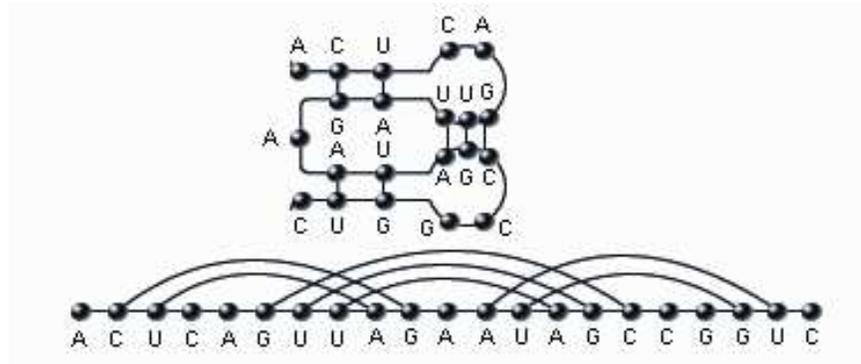,width=0.8\textwidth}\hskip15pt }
\caption{\small A $3$-noncrossing RNA structure with arc-length $\ge
3$, as a planar graph (top) and as a diagram (bottom)} \label{F:6}
\end{figure}

%%%
%%%%%%%%%%%%%%%%%%%%%%%%%%%%%%%%%%%%%%%%%%%%%%%%%%%%%%%%%%%%%%%%%%%%%%%%
%%%
Suppose we are given the primary RNA sequence
$$
{\bf A}{\bf C}{\bf U}{\bf C}{\bf A}{\bf G}{\bf U}{\bf U}{\bf A} {\bf
G}{\bf A}{\bf A}{\bf U}{\bf A}{\bf G}{\bf C}{\bf C}{\bf G}{\bf G}
{\bf U}{\bf C} \ .
$$
We then identify an RNA structure with the set of all bonds
different from the backbone-bonds of its primary sequence, i.e.~the
arcs $(i,i+1)$ for $1\le i\le n-1$. Accordingly an RNA structure is
a combinatorial graph over the labels of the nucleotides of the
primary sequence. These graphs can be represented in several ways.
In Figure~\ref{F:2} we represent a structure with loop-loop
interactions in two ways.
%%%%
%%%%%%%%%%%%%%%%%%%%%%%%%%%%%%%%%%%%%%%%%%%%%%%%%%%%%%%%%%%%%%%%%%%%%%%%%%
%%%%
%%%%
%%%%%%%%%%%%%%%%%%%%%%%%%%%%%%%%%%%%%%%%%%%%%%%%%%%%%%%%%%%%%%%%%%%%%%%%%%
%%%%
In the following we will consider structures as diagram
representations of digraphs. A digraph $D_n$ is a pair of sets
$V_{D_n},E_{D_n}$, where $V_{D_n}=
\{1,\dots,n\}$ and $E_{D_n}\subset \{(i,j)\mid 1\le i< j\le n\}$.
$V_{D_n}$ and $E_{D_n}$ are called vertex and arc set, respectively.
A $k$-noncrossing digraph is a digraph in which all
vertices have degree $\le 1$ and which does not contain a $k$-set of
arcs that are mutually intersecting, or equivalently
\begin{eqnarray}
\ \not\exists\,
(i_{r_1},j_{r_1}),(i_{r_2},j_{r_2}),\dots,(i_{r_k},j_{r_k});\quad & &
i_{r_1}<i_{r_2}<\dots<i_{r_k}<j_{r_1}<j_{r_2}<\dots<j_{r_k} \ .
\end{eqnarray}
We will represent digraphs as a diagrams (Figure~\ref{F:2}) by
representing the vertices as integers on a line and connecting any
two adjacent vertices by an arc in the upper-half plane. The
direction of the arcs is implicit in the linear ordering of the
vertices and accordingly omitted.
%%%
%%%%%%%%%%%%%%%%%%%%%%%%%%%%%%%%%%%%%%%%%%%%%%%%%%%%%%%%%%%%%%%%%%%%%%
%%%
\begin{definition}\label{D:rna}
An $k$-noncrossing RNA structure with arc-length $\ge 3$ is a digraph in
which all vertices have degree $\le 1$, having at most a $k-1$-set of
mutually intersecting arcs without arcs of length $\le 2$, i.e.~arcs of
the form $(i,i+1)$ and $(i,i+2)$, respectively.
Let ${\sf S}_{k,3}(n)$ and ${\sf S}_{k,3}(n,\ell)$ be the numbers of
$k$-noncrossing RNA structures with arc-length $\ge 3$ and those with
$\ell$ isolated vertices, respectively.
\end{definition}
%%%
%%%%%%%%%%%%%%%%%%%%%%%%%%%%%%%%%%%%%%%%%%%%%%%%%%%%%%%%%%%%%%%%%%%%%%
Let $f_{k}(n,\ell)$ denote the number of $k$-noncrossing digraphs with
$\ell$ isolated points. We have shown in \cite{Reidys:07pseu} that
\begin{align}\label{E:ww0}
f_{k}(n,\ell)& ={n \choose \ell} f_{k}(n-\ell,0) \\
\label{E:ww1}
\det[I_{i-j}(2x)-I_{i+j}(2x)]|_{i,j=1}^{k-1} &=
\sum_{n\ge 1} f_{k}(n,0)\,\frac{x^{n}}{n!} \\
\label{E:ww2}
e^{x}\det[I_{i-j}(2x)-I_{i+j}(2x)]|_{i,j=1}^{k-1}
&=(\sum_{\ell \ge 0}\frac{x^{\ell}}{\ell!})(\sum_{n \ge
1}f_{k}(n,0)\frac{x^{n}}{n!})=\sum_{n\ge 1}
\left\{\sum_{\ell=0}^nf_{k}(n,\ell)\right\}\,\frac{x^{n}}{n!} \ .
\end{align}
where $I_{r}(2x)=\sum_{j \ge 0}\frac{x^{2j+r}}{{j!(r+j)!}}$ is the
hyperbolic Bessel function of the first kind of order $r$. In
particular we have for $k=3$
\begin{equation}\label{E:2-3}
f_{3}(n,\ell)=
{n \choose \ell}\left[C_{\frac{n-\ell}{2}+2}C_{\frac{n-\ell}{2}}-
      C_{\frac{n-\ell}{2}+1}^{2}\right] \ ,
\end{equation}
where $C_m=\frac{1}{m+1}\binom{2m}{m}$ is the $m$th Catalan number.
The derivation of the generating function of $k$-noncrossing RNA structures,
given in Theorem~\ref{T:cool2} below uses advanced methods and novel
constructions of enumerative combinatorics due to Chen~{\it et.al.}
\cite{Chen:07a,Gessel:92a} and Stanley's mapping between matchings and
oscillating tableaux i.e.~families of Young diagrams in which any two
consecutive shapes differ by exactly one square.
The enumeration is obtained using the reflection principle due to Gessel
and Zeilberger \cite{Gessel:92a} and Lindstr\"om \cite{Lindstroem:73a}
in accord with an inclusion-exclusion argument in order to eliminate the
arcs of length $\le 2$ \cite{Reidys:07pseu}.
%%%
%%%%%%%%%%%%%%%%%%%%%%%%%%%%%%%%%%%%%%%%%%%%%%%%%%%%%%%%%%%%%%%%%%%%%%%%%%
%%%
\begin{theorem}\label{T:cool2}
Let $k\in\mathbb{N}$, $k>2$. Then the numbers of $k$-noncrossing
RNA structures ${\sf S}_{k,3}(n,\ell)$ and ${\sf S}_{k,3}(n)$ are given by
\begin{eqnarray}\label{E:da2}
{\sf S}_{k,3}(n,\ell) & = &
\sum_{b\ge 0}(-1)^{b}
\lambda(n,b)f_{k}(n-2b),\ell) \\
\label{E:da3}
{\sf S}_{k,3}(n) & = & \sum_{b=0}^{\lfloor n/2\rfloor}
(-1)^{b} \lambda(n,b)
\left\{\sum_{\ell=0}^{n-2b}f_{k}(n-2b,\ell)\right\}\ .
\end{eqnarray}
where $\lambda(n,b)$ satisfies the recursion
\begin{equation}\label{E:hh}
\lambda(n,b)=\lambda(n-1,b)+\lambda(n-2,b-1)+
\lambda(n-3,b-1)+\lambda(n-4,b-2)
\end{equation}
and the initial conditions for eq.~{\rm (\ref{E:hh})} are
$\lambda(n,0)=1$ and $\lambda(n,1)=2n-3$.
\end{theorem}
%%%
%%%%%%%%%%%%%%%%%%%%%%%%%%%%%%%%%%%%%%%%%%%%%%%%%%%%%%%%%%%%%%%%%%%%%%%%%%
%%%
%%%
%%%%%%%%%%%%%%%%%%%%%%%%%%%%%%%%%%%%%%%%%%%%%%%%%%%%%%%%%%%%%%%%%%%%%%%%%
%%%

\section{The exponential factor}\label{S:exp}

%%%
%%%%%%%%%%%%%%%%%%%%%%%%%%%%%%%%%%%%%%%%%%%%%%%%%%%%%%%%%%%%%%%%%%%%%%%%%
%%%
In this section we obtain the exponential growth factor of the
coefficients ${\sf S}_{k,3}(n)$. Let us begin by considering the
generating function $\sum_{n\ge 0}{\sf S}_{k,3}(n)x^n$ as a power series
over $\mathbb{R}$. Since $\sum_{n\ge 0}{\sf S}_k(n)x^n$ has
monotonously increasing coefficients $\lim_{n\to\infty}{\sf
S}_{k,3}(n)^{\frac{1}{n}}$ exists and determines via Hadamard's formula
its radius of convergence. Due to the inclusion-exclusion form of the
terms ${\sf S}_{k,3}(n)$, it is not obvious however, how to compute this
radius. Our strategy consists in first showing that ${\sf S}_{k,3}(n)$ is
closely related to $f_k(2n,0)$ via a functional relation of generating
functions.
%%%
%%%%%%%%%%%%%%%%%%%%%%%%%%%%%%%%%%%%%%%%%%%%%%%%%%%%%%%%%%%%%%%%%%%%%%%%%
%%%
\begin{lemma}\label{L:laplace}
Let $x$ be an indeterminate over $\mathbb{R}$ and $T_{k}(n)$ be the number of
$k$-noncrossing partial matchings over $[n]$, i.e. $T_{k}(n)=\sum_{m \le
\lfloor \frac{n}{2}\rfloor}\binom{n}{2m}f_k(2m,0)$.
Let furthermore $\rho_k$ denote the
radius of convergence of $\sum_{n \ge 0}T_{k}(n)\, x^n$.
Then we have
\begin{equation}\label{E:laplace}
\forall \, \vert x\vert <\rho_k;\quad
\sum_{n \ge 0}T_{k}(n)\, x^n =
\frac{1}{1-x}\, \sum_{n \ge 0}f_k(2n,0)\,\left(\frac{x}{1-x}\right)^{2n} \ .
\end{equation}
\end{lemma}
%%%
%%%%%%%%%%%%%%%%%%%%%%%%%%%%%%%%%%%%%%%%%%%%%%%%%%%%%%%%%%%%%%%%%%%%%%%%%
%%%
\begin{proof}
We will relate the ordinary generating functions of $T_{k}(n)$
\begin{equation}\label{E:ll}
\sum_{n \ge 0}T_{k}(n)\frac{x^n}{n!}=
\sum_{n \ge 0}\sum_{m\le\frac{n}{2}}{n \choose 2m}f_k(2m,0)
\frac{x^n}{n!}=e^x\cdot {\sf det}[I_{i-j}(2x)-I_{i+j}(2x)]_{i,j=1}^{k-1}\
\end{equation}
(see eq.~(\ref{E:ww2}))
and $f_k(2n,0)$ via Laplace transforms as follows: $\sum_{n \ge 0}
T_{k}(n)x^n = \sum_{n \ge 0}T_{k}(n)\frac{x^n}{n!}n!$ is convergent for
any $x\in\mathbb{R}$, whence we derive, using the Laplace transformation
and interchanging integration and summation
\begin{eqnarray*}
\sum_{n \ge 0}T_{k}(n)x^n  =  \sum_{n \ge 0}T_{k}(n)\frac{x^n}{n!}
                                \int_{0}^{\infty}e^{-t}t^n dt =
 \int_0^{\infty}\sum_{n \ge 0}T_{k}(n)\frac{(xt)^n}{n!}e^{-t}dt
\end{eqnarray*}
We next interpret the rhs via eq.~(\ref{E:ll}) and obtain
\begin{eqnarray*}
\sum_{n \ge 0}T_{k}(n)x^n
&=& \int_0^{\infty}e^{xt}{\sf det}[I_{i-j}(2xt)-I_{i+j}(2xt)]_{i,j=1}^{k-1}
    e^{-t}dt \ .
\end{eqnarray*}
Interchanging integration and summation yields
\begin{eqnarray*}
\sum_{n \ge 0}T_{k}(n)x^n
&=& \int_0^{\infty}e^{(x-1)t}\sum_{n \ge 0}f_k(2n,0)\frac{(xt)^{2n}}{(2n)!}dt\\
&=& \sum_{n \ge 0}f_k(2n,0)\cdot\frac{1}{(2n)!}
                  \int_0^{\infty}e^{(x-1)t}(xt)^{2n}dt\\
&=& \sum_{n \ge 0}f_k(2n,0)\frac{1}{(2n)!}\int_0^{\infty}
     e^{-(1-x)t}((1-x)t)^{2n}(\frac{x}{1-x})^{2n}dt\\
&=&\sum_{n \ge
0}f_k(2n,0)\frac{1}{(2n)!}\cdot(2n)!(\frac{x}{1-x})^{2n}\cdot\frac{1}{1-x}\\
&=&\frac{1}{1-x}\sum_{n \ge 0}f_k(2n,0)\,(\frac{x}{1-x})^{2n}
\end{eqnarray*}
and the proof of the lemma is complete.
\end{proof}
%%%
%%%%%%%%%%%%%%%%%%%%%%%%%%%%%%%%%%%%%%%%%%%%%%%%%%%%%%%%%%%%%%%%%%%%%%%%%
%%%
\begin{lemma}\label{L:func}
Let $x$ be an indeterminante over $\mathbb{R}$ and ${\sf S}_{k,3}(n)$
be the number of $k$-noncrossing RNA structures with arc-length $\ge
3$. Then we have the functional equation
\begin{equation}\label{E:laplace}
\sum_{n \ge 0}{\sf S}_{k,3}(n)\, x^n =
\frac{1}{1-x+x^2+x^3-x^4}\sum_{n \ge
0}f_k(2n,0)\left(\frac{x-x^3}{1-x+x^2+x^3-x^4}\right)^{2n}
\end{equation}
\end{lemma}
%%%
%%%%%%%%%%%%%%%%%%%%%%%%%%%%%%%%%%%%%%%%%%%%%%%%%%%%%%%%%%%%%%%%%%%%%%%%%
%%%
\begin{proof}
According to Theorem~\ref{T:cool2} we have
\begin{eqnarray*}
{\sf S}_{k,3}(n) & = & \sum_{b \le \lfloor \frac{n}{2}\rfloor}
(-1)^b\lambda(n,b)\sum_{\ell=0}^{n-2b}f_{k}(n-2b,\ell) \\
& = & \sum_{b \le \lfloor \frac{n}{2}\rfloor}
(-1)^b\lambda(n,b)\sum_{m=2b}^{n}
{n-2b \choose m-2b}f_{k}(m-2b,0)
\end{eqnarray*}
Our goal is now to relate $\sum_{n \ge 0}{\sf S}_{k,3}(n)x^n$ to the
terms $T_k(n)$. For this purpose we derive
\begin{align*}
\sum_{n \ge 0}{\sf S}_{k,3}(n)x^n&=\sum_{n \ge 0}\sum_{2b\le n}
(-1)^b\lambda(n,b)\sum_{m=2b}^{n}{n-2b \choose
m-2b}f_k(m-2b,0)\, x^n\\
&=\sum_{b \ge 0}(-1)^b x^{2b}\sum_{n \ge
2b}\lambda(n,b)T_{k}(n-2b)x^{n-2b}\\
&=\sum_{b \ge 0}(-1)^b x^{2b}\sum_{n \ge 0}\lambda(n+2b,b)T_{k}(n)\,x^n \ .\\
\end{align*}
Interchanging the summations w.r.t.~$b$ and $n$ we consequently arrive at
\begin{equation}\label{E:21}
\sum_{n \ge 0}{\sf S}_{k,3}(n)x^n = \sum_{n \ge 0}\left[\sum_{b\ge 0}(-1)^b
x^{2b}\lambda(n+2b,b)\right]T_{k}(n)\, x^n \ .
\end{equation}
We set $\varphi_{n}(x)=\sum_{b \ge 0}\lambda(n+2b,b)x^b$. According to
Theorem~\ref{T:cool2} we have the recursion formula
$$
\lambda(n+2b,b)=\lambda(n+2b-1,b)+\lambda(n+2b-2,b-1)+
\lambda(n+2b-3,b-1)+\lambda(n+2b-4,b-2)
$$
Multiplying with $x^b$ and taking the summation over all $b$
ranging from $0$ to $\lfloor n/2\rfloor$ implies the following functional
equation for $\varphi_{n}(x)$, $n=1,2\ldots$
\begin{equation}\label{E:22}
\varphi_{n}(x)=\varphi_{n-1}(x)+x\cdot \varphi_{n}(x)+x \cdot \varphi_{n-1}(x)+
x^2\varphi_{n}(x) \ .
\end{equation}
Eq.~(\ref{E:22}) is equivalent to
\begin{equation}\label{E:23}
\frac{\varphi_{n}(x)}{\varphi_{n-1}(x)}=\frac{1+x}{1-x-x^2}
\quad \text{\rm and }\quad \varphi_{0}(x)=\sum_{b \ge
0}\lambda(2b,b)x^b=\frac{1}{1-x-x^2}
\end{equation}
since $\lambda_{b}=\lambda(2b,b)$ satisfies the recursion formula
$\lambda_b=\lambda_{b-1}+\lambda_{b-2}$ and the initial condition
$\lambda_{0}=\lambda_{1}=1$. We can conclude from this that $\lambda(2b,b)$
is $b$-th Fibonacci number. As a result we obtain the formula
\begin{equation}\label{E:24}
\varphi_{n}(x)=\varphi_{0}(x)\left(\frac{1+x}{1-x-x^2}\right)^n=
\frac{1}{1-x-x^2}
\left(\frac{1+x}{1-x-x^2}\right)^n \ .
\end{equation}
Substituting eq.~(\ref{E:24}) into eq.~(\ref{E:21}) we can compute
\begin{eqnarray*}
\sum_{n \ge 0}{\sf S}_{k,3}(n)x^n &= & \sum_{n \ge
0}\varphi_{n}(-x^2)T_{k}(n)x^n \\
 & = &\sum_{n \ge
0}\frac{1}{1+x^2-x^4}\left(\frac{1-x^2}{1+x^2-x^4}\right)^n
T_{k}(n)\ x^n\\
& = & \frac{1}{1+x^2-x^4}
\sum_{n \ge
0} T_{k}(n)\, \left(\frac{x-x^3}{1+x^2-x^4}\right)^n \ .
\end{eqnarray*}
Via Lemma~\ref{L:laplace} we have the following interpretation of
$\sum_{n \ge 0}T_{k}(n) \, x^n$
$$
\sum_{n \ge 0}T_{k}(n) \, x^n = \frac{1}{1-x}\sum_{n \ge 0}
f_k(2n,0)\,(\frac{x}{1-x})^{2n}\ .
$$
Therefore we obtain setting $x'=\frac{x-x^3}{1+x^2-x^4}$, $\frac{1}{1-x'}=
\frac{1+x^2-x^4}{1-x+x^2+x^3-x^4}$ and
$\frac{x'}{1-x'}=\frac{x-x^3}{1-x+x^2+x^3-x^4}$
\begin{eqnarray*}
\sum_{n \ge 0}{\sf S}_{k,3}(n)x^n &= &
\frac{1}{1-x+x^2+x^3-x^4}\sum_{n \ge
0}f_k(2n,0)\left(\frac{x-x^3}{1-x+x^2+x^3-x^4}\right)^{2n} \ ,
\end{eqnarray*}
whence the lemma.
\end{proof}
%%%
%%%%%%%%%%%%%%%%%%%%%%%%%%%%%%%%%%%%%%%%%%%%%%%%%%%%%%%%%%%%%%%%%%%%%%%%%
%%%
Using complex analysis we can extend Lemma~\ref{L:func} to arbitrary
$\vert z\vert < \rho_k$, where $z\in \mathbb{C}$
%%%
%%%%%%%%%%%%%%%%%%%%%%%%%%%%%%%%%%%%%%%%%%%%%%%%%%%%%%%%%%%%%%%%%%%%%%%%%
%%%
\begin{lemma}\label{L:ana}
Let $k>2$ be an integer, then we have for arbitrary $z\in\mathbb{C}$
with the property $\vert z\vert <\rho_k$ the equality
\begin{equation}\label{E:rr2}
\sum_{n \ge 0}{\sf S}_{k,3}(n)\, z^n = \frac{1}{1-z+z^2+z^3-z^4}
\sum_{n \ge
0}f_k(2n,0)\left(\frac{z-z^3}{1-z+z^2+z^3-z^4}\right)^{2n} \ .
\end{equation}
\end{lemma}
%%%
%%%%%%%%%%%%%%%%%%%%%%%%%%%%%%%%%%%%%%%%%%%%%%%%%%%%%%%%%%%%%%%%%%%%%%%%%
%%%
\begin{proof}
The power series $\sum_{n\ge 0} {\sf S}_{k,3}(n) z^{n}$ and
$\frac{1}{1-z+z^2+z^3-z^4}
\sum_{n\ge 0} f_k(2n,0) \left(\frac{z-z^3}{1-z+z^2+z^3-z^4}\right)^{2n}$
are analytic in a disc of radius $0<\epsilon<\rho_k$ and
according to Lemma~\ref{L:func} coincide on the interval
$]-\epsilon,\epsilon [$. Therefore both functions are analytic and equal
on the sequence $(\frac{1}{n})_{n\in\mathbb{N}}$
which converges to $0$ and
standard results of complex analysis (zeros of nontrivial
analytic functions are isolated) imply that eq.~(\ref{E:rr2})
holds for any $z\in\mathbb{C}$ with $\vert z\vert<\rho_k$, whence
the lemma.
\end{proof}

Lemma~\ref{L:ana} is the key to compute the exponential growth rates
for any $k>2$. In its proof we recruit the Theorem of Pringsheim
\cite{Titmarsh:39} which asserts that a power series $\sum_{n\ge
0}a_nz^n$ with $a_n\ge 0$ has its radius of convergence as dominant
(but not necessarily unique) singularity. In particular there exists
a dominant real valued singularity.
%%%
%%%%%%%%%%%%%%%%%%%%%%%%%%%%%%%%%%%%%%%%%%%%%%%%%%%%%%%%%%%%%%%%%%%%%%%%%
%%%
\begin{theorem}\label{T:asy1}
Let $k\ge 3$ be a positive integer and $r_k$ be the radius of
convergence of the power series $\sum_{n\ge 0}f_k(2n,0)z^{2n}$ and
\begin{equation}\label{E:theta}
\vartheta\colon [0, \frac{\sqrt{2}}{2}]\longrightarrow
[0,\frac{5-\sqrt{2}}{4}], \quad z\mapsto
\frac{z(1-z)(1+z)}{-(z^2-\frac{1}{2})^2+z(z^2-\frac{1}{2})-
\frac{z}{2}+\frac{5}{4}} \ .
\end{equation}
Then the power series $\sum_{n\ge 0}{\sf S}_{k,3}(n)z^n$ has the real valued,
dominant singularity $\rho_k$, which is the unique real solution
of $\vartheta(x)=r_k$ and for the number of
$k$-noncrossing RNA structures with arc-length $\ge 3$ holds
\begin{equation}\label{E:rel}
{\sf S}_{k,3}(n)\sim \left(\frac{1}{\rho_k}\right)^n \ .
\end{equation}
\end{theorem}
%%%
%%%%%%%%%%%%%%%%%%%%%%%%%%%%%%%%%%%%%%%%%%%%%%%%%%%%%%%%%%%%%%%%%%%%%%%%%
%%%
In Section~\ref{S:sub} we will in particular prove that
$\rho_3\approx 0.21982$.
\begin{proof}
Suppose we are given $r_k$, then $r_k\le \frac{1}{2}$ (this follows
immediately from $C_n\sim 2^{2n}$ via Stirling's formula). The
functional identity of Lemma~\ref{L:func} allows us to derive the
radius of convergence of $\sum_{n\ge 0}S_k(n)z^n$. According to
Lemma~\ref{L:ana} we have
\begin{equation}\label{E:25}
\sum_{n \ge 0}{\sf S}_{k,3}(n)\, z^n = \frac{1}{1-z+z^2+z^3-z^4}
\sum_{n \ge
0}f_k(2n,0)\left(\frac{z-z^3}{1-z+z^2+z^3-z^4}\right)^{2n}
\end{equation}
$f_k(2n,0)$ is monotone, whence the limit $\lim_{n\to
\infty}f_k(2n,0)^{ \frac{1}{2n}}$ exists and applying Hadamard's
formula: $\lim_{n\to
\infty}f_k(2n,0)^{\frac{1}{2n}}=\frac{1}{{r_{k}}}$. For $z\in
\mathbb{R}$, we proceed by computing the roots of
$$
\left|\frac{z-z^3}{1-z+z^2+z^3-z^4}\right|={r_{k}}
$$
which for $r_k\le \frac{1}{2}$ has the minimal root $\rho_k$. We
next show that $\rho_k$ is indeed the radius of convergence of
$\sum_{n\ge 0} {\sf S}_k(n) z^n$. For this purpose we observe that
the map
\begin{equation}\label{E:w=1}
\vartheta\colon [0, \frac{\sqrt{2}}{2}]\longrightarrow
[0,\frac{5-\sqrt{2}}{4}], \quad z\mapsto
\frac{z(1-z)(1+z)}{-(z^2-\frac{1}{2})^2+z(z^2-\frac{1}{2})-
\frac{z}{2}+\frac{5}{4}}
, \qquad \text{\rm where} \quad\vartheta(\rho_k)={r_k}
\end{equation}

\begin{figure}[ht]
\centerline{ \epsfig{file=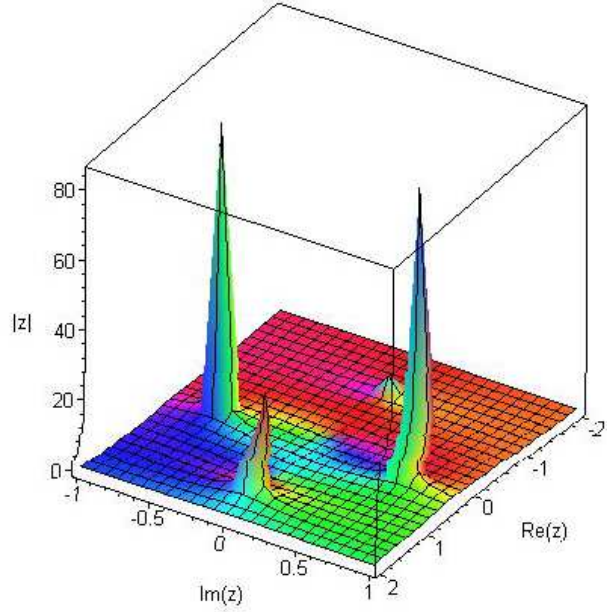,width=0.6\textwidth}\hskip15pt }
\caption{\small We display 4 poles(the corresponding 4 peaks) in
$\frac{z-z^3}{1-z+z^2+z^3-z^4}$ over $\mathbb{C}$. The picture
illustrates that $\vartheta$ is a bijection on the interval $[0,
\frac{\sqrt{2}}{2}]$, which allows us to obtain the dominant singularity
$\rho_k$.} \label{F:6}
\end{figure}

is bijective, continuous and strictly decreasing. Continuity and
strict monotonicity of $\vartheta$ guarantee in view of
eq.~(\ref{E:w=1}) that $\rho_k$, is indeed the radius of convergence
of the power series $\sum_{n\ge 0} {\sf S}_k(n) z^n$. In order to
show that $\rho_k$ is a dominant singularity we consider $\sum_{n\ge
0}{\sf S}_{k,3}(n)z^n$ as a power series over $\mathbb{C}$. Since
${\sf S}_{k,3}(n)\ge 0$, the theorem of
Pringsheim~\cite{Titmarsh:39} guarantees that $\rho_k$ itself is a
singularity. By construction $\rho_k$ has minimal absolute value and
is accordingly dominant. Since ${\sf S}_{k,3}(n)$ is monotone
$\lim_{n\to\infty}{\sf S}_{k,3}(n)^{\frac{1} {n}}$ exists and we
obtain using Hadamard's formula
\begin{equation}
\lim_{n\to\infty}{\sf
S}_{k,3}(n)^{\frac{1}{n}}=\frac{1}{\rho_k},\quad \text{\rm or
equivalently}\quad {\sf S}_{k,3}(n)\sim \left(\frac{1}{\rho_k}
\right)^n \, ,
\end{equation}
from which eq.~(\ref{E:rel}) follows and the proof of the theorem is complete.
\end{proof}
%%%%
%%%%%%%%%%%%%%%%%%%%%%%%%%%%%%%%%%%%%%%%%%%%%%%%%%%%%%%%%%%%%%%%%%%%%%%%%
%%%
\section{Asymptotic Analysis}\label{S:sub}
%%%
%%%%%%%%%%%%%%%%%%%%%%%%%%%%%%%%%%%%%%%%%%%%%%%%%%%%%%%%%%%%%%%%%%%%%%%%%
%%%

\begin{figure}[ht]
\centerline{%
\epsfig{file=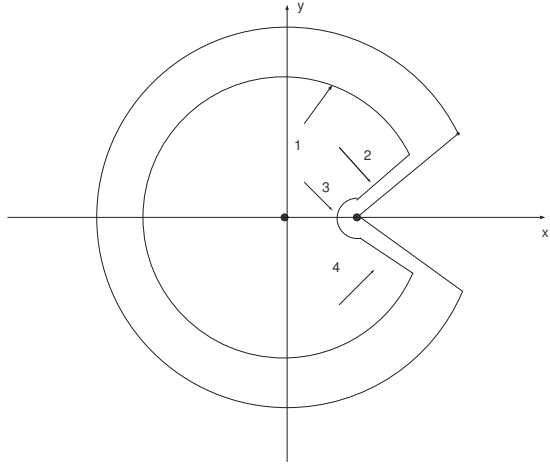,width=0.5\textwidth}\hskip15pt } \caption{\small
$\Delta_1$-domain enclosing a Hankel contour. We assume $z=1$ to be
the unique dominant singularity. The coefficients are obtained via
Cauchy's integral formula and the integral path is decomposed in $4$
segments. Segment $1$ becomes asymptotically irrelevant since by
construction the function involved is bounded on this segment.
Relevant are the rectilinear segments $2$ and $4$ and the inner
circle $3$. The only contributions to the contour integral are being
made here, which shows why the singular expansion allows to
approximate the coefficients so well.} \label{F:7}
\end{figure}

In this section we provide the asymptotic number of $3$-noncrossing RNA
structures with arc-length $\ge 3$.
In the course of our analysis we derive the analytic continuation of the
power series $\sum_{n\ge 0}{\sf S}_{3,3}(n)z^n$. The analysis will in
particular provide independent proof of the exponential factor computed
in Theorem~\ref{T:asy1}. The derivation of the subexponential factors is
based on singular expansions \cite{Flajolet:05} in combination with transfer
theorems. The key ingredient for the coefficient extraction is the Hankel
contours, see Figure~\ref{F:7}. Let us begin by specifying a
suitable domain for our Hankel contours tailored for
Theorem~\ref{T:transfer1}.

\begin{definition}\label{D:delta}
Given two numbers $\phi,R$, where $R>1$ and $0<\phi<\frac{\pi}{2}$
and $\rho\in\mathbb{R}$ the open domain $\Delta_\rho(\phi,R)$ is
defined as
\begin{equation}
\Delta_\rho(\phi,R)=\{ z\mid \vert z\vert < R, z\neq \rho,\, \vert
{\rm Arg}(z-\rho)\vert >\phi\}
\end{equation}
A domain is a $\Delta_\rho$-domain if it is of the form
$\Delta_\rho(\phi,R)$ for some $R$ and $\phi$. A function is
$\Delta_\rho$-analytic if it is analytic in some
$\Delta_\rho$-domain.
\end{definition}
%%%
%%%%%%%%%%%%%%%%%%%%%%%%%%%%%%%%%%%%%%%%%%%%%%%%%%%%%%%%%%%%%%%%%%%%%%%%%
%%%
%%%%
%%%%%%%%%%%%%%%%%%%%%%%%%%%%%%%%%%%%%%%%%%%%%%%%%%%%%%%%%%%%%%%%%%%%%%%%%%
%%%%
%%%%
%%%%%%%%%%%%%%%%%%%%%%%%%%%%%%%%%%%%%%%%%%%%%%%%%%%%%%%%%%%%%%%%%%%%%%%%%%
%%%%
Since the Taylor coefficients have the property
\begin{equation}\label{E:scaling}
\forall \,\gamma\in\mathbb{C}\setminus 0;\quad [z^n]f(z)=\gamma^n
[z^n]f(\frac{z}{\gamma}) \ ,
\end{equation}
we can, w.l.o.g.~reduce our analysis to the case where $1$ is the
dominant singularity. We use $U(a,r)=\{z\in \mathbb{C}|\vert
z-a\vert<r\}$ to denote the open neighborhood of $a$ in
$\mathbb{C}$.

 We use the notation
\begin{equation}\label{E:genau}
\left(f(z)=O\left(g(z)\right) \ \text{\rm as $z\rightarrow
\rho$}\right)\quad \Longleftrightarrow \quad \left(f(z)/g(z) \
\text{\rm is bounded as $z\rightarrow \rho$}\right)
\end{equation}
and if we write $f(z)=O(g(z))$ it is implicitly assumed that $z$
tends to a (unique) singularity. $[z^n]\,f(z)$ denotes the
coefficient of $z^n$ in the power series expansion of $f(z)$ around
$0$.
%%%
%%%%%%%%%%%%%%%%%%%%%%%%%%%%%%%%%%%%%%%%%%%%%%%%%%%%%%%%%%%%%%%%%%%%%%%%%
%%%
\begin{theorem}\label{T:transfer1}\cite{Flajolet:05}
Let $\alpha$ be an arbitrary complex number in $\mathbb{C}\setminus
\mathbb{Z}_{\le 0}$, $f(z)$ be a $\Delta_1$-analytic function.
Suppose $r\in\mathbb{Z}_{\ge 0}$, and
$f(z)=O((1-z)^{r}\ln^{}(\frac{1}{1-z}))$ in the intersection of a
neighborhood of $1$ and the $\Delta_1$-domain, then we have
\begin{equation}
[z^n]f(z)\sim K\,  (-1)^r\frac{r!}{n(n-1)\dots(n-r)} \quad \text{\it
for some $K>0$}\ .
\end{equation}
\end{theorem}
%%%
%%%%%%%%%%%%%%%%%%%%%%%%%%%%%%%%%%%%%%%%%%%%%%%%%%%%%%%%%%%%%%%%%%%%%
%%%

We are now prepared to compute an explicit formula for the numbers
of $3$-noncrossing RNA structures with arc-length $\ge 3$.

%%%
%%%%%%%%%%%%%%%%%%%%%%%%%%%%%%%%%%%%%%%%%%%%%%%%%%%%%%%%%%%%%%%%%%%%%%%%%
%%%
\begin{theorem}\label{T:asy3}
The number of $3$-noncrossing RNA structures with arc length $\ge 3$
is asymptotically given by
\begin{eqnarray*}
\label{E:konk3} {\sf S}_{3,3}(n) & \sim & \frac{6.11170\cdot
4!}{n(n-1)\dots(n-4)}\,
4.54920^n \ .\\
\end{eqnarray*}
\end{theorem}
%%%
%%%%%%%%%%%%%%%%%%%%%%%%%%%%%%%%%%%%%%%%%%%%%%%%%%%%%%%%%%%%%%%%%%%%%%%%%
%%%
\begin{proof}
{\it Claim $1$.}
The dominant singularity $\rho_3$ of the power series
$\sum_{n\ge 0} {\sf S}_{3,3}(n) z^n$ is unique.\\
%%%%%%%%%%%%%%%%%%%%%%%%%%%%%%%%%%%%%%%%%%%
In order to prove Claim $1$ we use Lemma~\ref{L:ana}, according to
which the analytic function $\Xi_3(z)$ is the analytic continuation
of the power series $\sum_{n\ge 0} {\sf S}_{3}(n) z^n$. We proceed
by showing that $\Xi_3(z)$ has exactly $12$ singularities in
$\mathbb{C}$ and the dominant singularity is unique. The first four
singularities are the roots of the quartic polynomial
$P(z)=1-z+z^2+z^3-z^4$. Next we observe in analogy to our proof
in \cite{Reidys:07asy} that the power series
$\sum_{n\ge 0} f_3(2n,0) y^{n}$ has the analytic continuation
$\Psi(y)$ (obtained by MAPLE sumtools) given by
\begin{equation}\label{E:psi}
\Psi(y)=
\frac{-(1-16y)^{\frac{3}{2}}  P_{3/2}^{-1}(-\frac{16y+1}{16y-1})}
{16\, {y}^{\frac{5}{2}}} \ ,
\end{equation}
where $P_{\nu}^{m}(x)$ denotes the Legendre Polynomial of the first
kind with the parameters $\nu=\frac{3}{2}$ and $m=-1$. $\Psi(y)$ has
one dominant singularity at $y=\frac{1}{16}$, which in view of
$\vartheta(z)=(\frac{z-z^3}{1+z^2-z^4-z+z^3})^2$ induces exactly $8$
singularities of
$\Xi_3(z)=\frac{1}{1+z^2-z^4-z+z^3}\,
   \Psi\left(\left(\frac{z-z^3}{1+z^2-z^4-z+z^3}\right)^2\right)
$. Indeed, $\Psi(y^2)$ has the two singularities $\mathbb{C}$:
$\beta_1=\frac{1}{4}$ and $\beta_2=-\frac{1}{4}$ which produces
for $\Xi_3(z)$ (solving the quartic equation) the 8 singularities
$\rho_3\approx 0.21982$, $\zeta_2\approx
5.00829$, $\zeta_3\approx -1.07392$, $\zeta_4\approx 0.84581$,
$\zeta_5\approx -0.53243+0.11951i$, $\zeta_6 \approx -0.53243-0.11951i$,
$\zeta_7 \approx 1.10477$ and $\zeta_8 \approx -3.03992$.
The above values have error terms
(for details for solving the general quartic equation see Section~\ref{S:app})
of the order $10^{-5}$, which allows us
to conclude that the dominant singularity $\rho_3$ is unique and Claim
$1$ follows. \\
%%%%%%%%%%%%%%%%%%%%%%%%%%%%%%%%%%%%%%%%%%%%%%%%%%%%%%%%%%%%%%%%%%%%%%%%%
{\it Claim $2$.} \cite{Reidys:07asy} $\Psi(z)$ is
$\Delta_{\frac{1}{16}}(\phi,R)$-analytic and has the singular
expansion $(1-16z)^4\ln\left(\frac{1}{1-16z}\right)$.
\begin{equation}
\forall\, z\in\Delta_{\frac{1}{16}}(\phi,R)\cap
U(\frac{1}{16},\epsilon);\quad
\Psi(z)={O}\left((1-16z)^4\ln\left(\frac{1}{1-16z}\right)\right) \ .
\end{equation}
First $\Delta_{\frac{1}{16}}(\phi,R)$-analyticity of the function
$\Psi(z)$ is obvious.
We proceed by proving that $(1-16z)^4\ln\left(\frac{1}{1-16z}\right)$ is
its singular expansion in the intersection of a neighborhood of $\frac{1}{16}$
and the $\Delta$-domain $\Delta_{\frac{1}{16}}(\phi,R)$.
Using the notation of falling factorials
$(n-1)_4=(n-1)(n-2)(n-3)(n-4)$ we observe
$$
f_3(2n,0)=C_{n+2}C_{n}-C_{n+1}^2= \frac{1}{(n-1)_4}
\frac{12(n-1)_4(2n+1)}{(n+3)(n+1)^2(n+2)^2}\,
\binom{2n}{n}^2 \ .
$$
With this expression for $f_3(2n,0)$ we arrive at the formal identity
\begin{eqnarray*}
\sum_{n\ge 5}16^{-n}f_3(2n,0)z^n  & = &
O(\sum_{n\ge 5}
\left[16^{-n}\,\frac{1}{(n-1)_4}
\frac{12(n-1)_4(2n+1)}{(n+3)(n+1)^2(n+2)^2}\,
\binom{2n}{n}^2-\frac{4!}{(n-1)_4}\frac{1}{\pi}\frac{1}{n}\right]z^n \\
& & + \sum_{n\ge 5}\frac{4!}{(n-1)_4}\frac{1}{\pi}\frac{1}{n}z^n) \ ,
\end{eqnarray*}
where $f(z)=O(g(z))$ denotes that the limit $f(z)/g(z)$ is bounded
for $z\rightarrow 1$, eq.~(\ref{E:genau}). It is clear that
\begin{eqnarray*}
& & \lim_{z\to 1}(\sum_{n\ge 5}\left[16^{-n}\,\frac{1}{(n-1)_4}
\frac{12(n-1)_4(2n+1)}{(n+3)(n+1)^2(n+2)^2}\,
\binom{2n}{n}^2-\frac{4!}{(n-1)_4}\frac{1}{\pi}\frac{1}{n}\right]z^n)  \\
&= &
\sum_{n\ge 5} \left[16^{-n}\,\frac{1}{(n-1)_4}
\frac{12(n-1)_4(2n+1)}{(n+3)(n+1)^2(n+2)^2}\,
\binom{2n}{n}^2-\frac{4!}{(n-1)_4}\frac{1}{\pi}\frac{1}{n}\right]
 <\kappa
\end{eqnarray*}
for some $\kappa< 0.0784$. Therefore we can conclude
\begin{equation}
\sum_{n\ge 5}16^{-n}f_3(2n,0)z^n=
O(\sum_{n\ge 5}\frac{4!}{(n-1)_4}\frac{1}{\pi}\frac{1}{n}z^n) \ .
\end{equation}
We proceed by interpreting the power series on the rhs, observing
\begin{equation}
\forall\, n\ge 5\, ; \qquad
[z^n]\left((1-z)^4\,\ln\frac{1}{1-z}\right)=
\frac{4!}{(n-1)\dots (n-4)}\frac{1}{n} \, ,
\end{equation}
whence $\left((1-z)^4\,\ln\frac{1}{1-z}\right)$ is
the unique analytic continuation of $\sum_{n\ge 5}\frac{4!}{(n-1)_4}
\frac{1}{\pi}\frac{1}{n}z^n$.
Using the scaling property of Taylor coefficients
$[z^n]f(z)=\gamma^n [z^n]f(\frac{z}{\gamma})$ we obtain
\begin{equation}\label{E:isses}
\forall\, z\in\Delta_{\frac{1}{16}}(\phi,R)\cap
U(\frac{1}{16},\epsilon);\quad \Psi(z)
=O\left((1-16z)^4\ln\left(\frac{1}{1-16z}\right)\right) \ .
\end{equation}
Therefore we have proved that $(1-16z)^{4}\ln^{}(\frac{1}{1-16z})$
is the singular expansion of $\Psi(z)$ at $z=\frac{1}{16}$, whence
Claim $2$. Our last step consists in verifying that the type of the
singularity does not change when passing from $\Psi(z)$ to
$\Xi_3(z)= \frac{1}{1-z+z^2+z^3-z^4}
\Psi((\frac{z-z^3}{1-z+z^2+z^3-z^4})^2)$.
\\
%%%%%%%%%%%%%%%%%%%%%%%%%%%%%%%%%%%%%%%%%%%%%%%%%%%%%%%%%%%%%%%%%%%%%%%%%%%
{\it Claim $3$.} For $z\in
\Delta_{\rho_3}(\phi,R)\cap U(\rho_3,\epsilon)$ we have $\Xi_3(z)
={O}\left((1-\frac{z}{\rho_3})^4\ln(\frac{1}{1-\frac{z}{\rho_3}})
\right)$.\\
To prove the claim set $u(z)=1-z+z^2+z^3-z^4$. We first observe that
Claim $2$ and Lemma~\ref{L:ana} imply
\begin{align*}
\Xi_3(z) &=O\left(
\frac{1}{u(z)}\,
\left[\left(1-16(\frac{z-z^3}{u(z)})^2\right)^4
\ln\frac{1}{\left(1-16(\frac{z-z^3}{u(z)})^2\right)}\right]\right)
\ .
\end{align*}
The Taylor expansion of $q(z)=1-16(\frac{z-z^3}{u(z)})^2$ at
$\rho_3$ is given by $q(z)=\alpha(\rho_3-z)+{O}(z-\rho_3)^2$ and
setting $\alpha\approx -1.15861$ we compute
\begin{align*}
\frac{1}{u(z)}\, \left[q(z)^4\ln\frac{1} {q(z)}\right] &=
\frac{(\alpha(\rho_3-z)+{O}(z-\rho_3)^2)^4\ln\frac{1}{\alpha(\rho_3-z)+{O}
(z-\rho_3)^2}}{0.83679-0.45789(z-\rho_3)+O((z-\rho_3)^2)}\\
&= \frac{\left([\alpha+O(z-\rho_3)](\rho_3-z)^4
\ln\frac{1}{[\alpha+O(z-\rho_3)](\rho_3-z)} \right)}{O(z-\rho_3)}
\\
&={O}((\rho_3-z)^4\ln\frac{1}{\rho_3-z}) \ ,
\end{align*}
whence Claim $3$. Now we are in the position to employ
Theorem~\ref{T:transfer1}, and obtain for ${\sf S}_{3,3}(n)$
\begin{align*}
{\sf S}_{3,3}(n)&\sim K'\,
[z^n]\left((\rho_3-z)^4\ln\frac{1}{\rho_3-z} \right) \sim K'\,
\frac{4!} {n(n-1)\dots(n-4)}\left(\frac{1}{\rho_3}\right)^n  \ .
\end{align*}
Theorem~\ref{T:cool2} allows us to compute $K'=6.11170$ and the proof of
the Theorem is complete.
\end{proof}
%%%
%%%%%%%%%%%%%%%%%%%%%%%%%%%%%%%%%%%%%%%%%%%%%%%%%%%%%%%%%%%%%%%%%%%%%
%%%

\section{Appendix}\label{S:app}

%%%
%%%%%%%%%%%%%%%%%%%%%%%%%%%%%%%%%%%%%%%%%%%%%%%%%%%%%%%%%%%%%%%%%%%%%
%%%

Let us first introduce some basic definitions for quartic equation
(quartic) used in the following. $Ax^4+Bx^3+Cx^2+Dx+E=0$ is
a quartic if $A \ne 0$. A depressed quartic is a quartic such that
$B=0$ holds. $Ax^3+Bx^2+Cx+d=0$ is a cubic equation if $A \ne 0$ and in
particular a cubic without $x^2$ term is called as depressed cubic.
\\
The first step in solving the quartic consists in transforming it into
a depressed quartic. i.e. eliminate the $x^3$ term. Transform the
equation into $x^4+\frac{B}{A}x^3+\frac{C}{A}x^2+\frac{D}{A}x+\frac{E}{A}=0$
and substitute $x=u-\frac{B}{4A}$. Simplifying the original quartic yields
\begin{equation}\label{A:depqua}
u^4+\alpha u^2+\beta u +\gamma=0
\end{equation}
where $\alpha=\frac{-3B^2}{8A^2}+\frac{C}{A}$,
$\beta=\frac{B^3}{8A^3}-\frac{BC}{2A^2}+\frac{D}{A}$ and
$\gamma=\frac{-3B^4}{256A^4}+\frac{CB^2}{16A^3}-\frac{BD}{4A^2}
+\frac{E}{A}$, and eq.~(\ref{A:depqua}) is a depressed quartic
function, which is tantamount to
\begin{equation}\label{A:preu}
(u^2+\alpha)^2+\beta u+\gamma=\alpha u^2+\alpha^2
\end{equation}
The next step is to insert a variable $y$ into the perfect square
on the left side of eq.~(\ref{A:preu}). Add both
\begin{align*}
(u^2+\alpha+y)^2-(u^2+\alpha)^2=2yu^2+2y\alpha+y^2\\
0=(\alpha+2y)u^2-2yu^2-\alpha u^2
\end{align*}
to eq.(\ref{A:preu}) yields
\begin{equation}\label{A:insert-y}
(u^2+\alpha+y)^2=(\alpha+2y)u^2-\beta
u+(y^2+2y\alpha+\alpha^2-\gamma) \ .
\end{equation}
The next step consists in substituting for $y$ such that the right
side of eq.~(\ref{A:insert-y}) becomes a square. Observe that
$(su+t)^2=(s^2)u^2+(2st)u+(t^2)$ holds for any $s$ and $t$, and
the relation between the coefficients of the rhs is
$(2st)^2=4(s^2)(t^2)$. Therefore to make the rhs of
eq.~(\ref{A:insert-y}) into a perfect square, the following
equation must hold.
\begin{equation}{\label{A:solve-y}}
2y^3+5 \alpha y^2+(4 \alpha^2-2 \gamma)y+(\alpha^3-\alpha
\gamma-\frac{\beta^2}{4})=0 \ .
\end{equation}
Similarly, transform eq~(\ref{A:solve-y}) into a
depressed cubic equation by substituting $y=v-\frac{5}{6}\alpha$
\begin{equation}\label{A:depcub}
v^3+(-\frac{\alpha^2}{12}-\gamma)v+(-\frac{\alpha^3}{108}+\frac{\alpha
\gamma}{3}-\frac{\beta^2}{8})=0 \ .
\end{equation}
Set $P=-\frac{\alpha^2}{12}-\gamma$
and $Q=-\frac{\alpha^3}{108}+\frac{\alpha
\gamma}{3}-\frac{\beta^2}{8}$. Select any of the solutions of
eq.~(\ref{A:depcub}) of the form $v=\frac{P}{3U}-U$ where
$U=\sqrt[3]{\frac{Q}{2}\pm\sqrt{\frac{Q^2}{4}+\frac{P^3}{27}}}$
and $U \ne 0$, otherwise $v=0$. In view of $y=v-\frac{5}{6}$ this solution
yields for eq.~(\ref{A:solve-y})
$y=-\frac{5}{6}\alpha+\frac{P}{3U}-U$ for $U \ne 0$ and
$y=-\frac{5}{6}\alpha$ for $U=0$. Now the rhs of
eq.~(\ref{A:insert-y}) becomes
\begin{align*}
(\alpha+2y)u^2+(-\beta)u+(y^2+2y\alpha+\alpha^2-\gamma)=\left(
\left(\sqrt{\alpha+2y}\right)u+\frac{-\beta}{2\sqrt{\alpha+2y}}\right)^2
\end{align*}
Combined with eq.~(\ref{A:insert-y}) this yields two solutions for $u$:
$$
u=\frac{\pm\sqrt{\alpha+2y}\pm\sqrt{-\left(3\alpha+2y\pm\frac{2\beta}
{\sqrt{\alpha+2y}}\right)}}{2}
$$
where the first and the third $\pm$ must have the same sign.
This allows us to obtain the solutions for $x$:
$$
x=-\frac{B}{4A}+\frac{\pm\sqrt{\alpha+2y}\pm\sqrt{-\left(3\alpha+
                           2y\pm\frac{2\beta}{\sqrt{\alpha+2y}}\right)}}{2} \ .
$$
In particular for $x^4-5x^3-x^2+5x-1=0$, $A=1,B=-5,C=-1,D=5,E=-1$,
and hence $\alpha=-\frac{83}{8}$, $\beta=-\frac{105}{8}$,
$P=-\frac{16}{3}$ and $Q=\frac{299}{216}$ and $U \approx
1.21481-0.54955i \ne 0$, therefore $y \approx 6.21621$, and the
solutions are $\rho_3\approx 0.21982$, $\zeta_2\approx 5.00829$,
$\zeta_3\approx -1.07392$, $\zeta_4\approx 0.84581$. As for the
equation $x^4+3x^3-x^2-3x-1=0$, the corresponding solutions are
$\zeta_5\approx -0.53243+0.11951i$, $\zeta_6 \approx
-0.53243-0.11951i$, $\zeta_7 \approx 1.10477$ and $\zeta_8 \approx
-3.03992$.

%%%
%%%%%%%%%%%%%%%%%%%%%%%%%%%%%%%%%%%%%%%%%%%%%%%%%%%%%%%%%%%%%%%%%%%%%%%%%%
%%%
{\bf Acknowledgments.}
%%%
%%%%%%%%%%%%%%%%%%%%%%%%%%%%%%%%%%%%%%%%%%%%%%%%%%%%%%%%%%%%%%%%%%%%%%%%%%
%%%
%We are grateful to ??? for helpful comments. 
This work was supported by the
973 Project, the PCSIRT Project of the Ministry of Education, the Ministry
of Science and Technology, and the National Science Foundation of China.

\bibliography{a2}
\bibliographystyle{plain}

%%%
%%%%%%%%%%%%%%%%%%%%%%%%%%%%%%%%%%%%%%%%%%%%%%%%%%%%%%%%%%%%%%%%%%%%%%%%%%
%%%

\end{document}